%% file: Vrettos_et_al_TSG_2016_p1_arXiv.tex
\documentclass[11pt,a4paper]{article}

\usepackage{cite}
\usepackage[cmex10]{amsmath}
\usepackage{amsmath}
\usepackage{amssymb}
\usepackage{booktabs}
\usepackage{color}
\usepackage{times}
\usepackage{paralist}
\usepackage{subfigure}
\usepackage{float}
\usepackage{epstopdf}
\usepackage{authblk}
\usepackage{graphicx}

\let\mathnumsetfont\mathbb
\newcommand\Nset{\mathnumsetfont N}       

\usepackage[utf8]{inputenc}
\usepackage[cmex10]{amsmath}
\usepackage{amssymb}
\usepackage{mathtools}
\usepackage{subfigure}
\usepackage{float}
\usepackage{booktabs}
\usepackage{times}
\usepackage{paralist}
\usepackage{cite}
\usepackage{algorithm}
\usepackage{algpseudocode}
\usepackage{epstopdf}
\usepackage{soul} 
\usepackage{color}
\usepackage{eurosym}
\usepackage[printonlyused]{acronym}
\graphicspath{
    {FiguresBldsExp/}
}
\DeclareGraphicsExtensions{.pdf,.png,.jpg,.jpeg,.eps}

\usepackage{amsthm}

\newtheorem{proposition}{Proposition}
\newtheorem{assumption}{Assumption}
\newtheorem{lemma}{Lemma}

\definecolor{Red}{rgb}{1,0,0}
\definecolor{Green}{rgb}{0,1,0}
\definecolor{Blue}{rgb}{0,0,1}

\begin{document}
\title{Experimental Demonstration of Frequency Regulation by Commercial Buildings -- Part I: Modeling and Hierarchical Control Design}
\author[1]{Evangelos~Vrettos\thanks{vrettos@eeh.ee.ethz.ch}}
\author[2]{Emre~C.~Kara\thanks{eckara@slac.stanford.edu}}
\author[3]{Jason~MacDonald\thanks{jsmacdonald@lbl.gov}}
\author[1]{G\"{o}ran~Andersson\thanks{andersson@eeh.ee.ethz.ch}}
\author[4]{Duncan~S.~Callaway\thanks{dcal@berkeley.edu}}
\affil[1]{Power Systems Laboratory, ETH Zurich, Switzerland}
\affil[2]{National Accelerator Laboratory (SLAC), California, US}
\affil[3]{Grid Integration Group, Lawrence Berkeley National Laboratory (LBNL), California, US}
\affil[4]{Energy and Resources Group, University of California, Berkeley, US}
\maketitle

\begin{abstract}
This paper is the first part of a two-part series in which we present results from an experimental demonstration of frequency regulation in a commercial building test facility.

In Part I, we introduce the test facility and develop relevant building models. Furthermore, we design a hierarchical controller that consists of three levels: a reserve scheduler, a building climate controller, and a fan speed controller for frequency regulation. We formulate the reserve scheduler as a robust optimization problem and introduce several approximations to reduce its complexity. The building climate controller is comprised of a robust model predictive controller and a Kalman filter. The frequency regulation controller consists of a feedback and a feedforward loop, provides fast responses, and is stable.

Part I presents building model identification and controller tuning results, whereas Part II reports results from the operation of the hierarchical controller under frequency regulation.
\end{abstract}

\newpage
\input{acronymDefinitions_sorted.tex}

\input{introduction_p1.tex}
\input{expFacility_v3.tex}
\input{modeling_v3.tex}
\input{controlDesign_v3.tex}
\input{conclusion_p1.tex}

\begin{appendix}
\input{appendix_proofs.tex}
\end{appendix}

\bibliographystyle{elsarticle-num}
\bibliography{biblio_EV}

\end{document}

%% file: acronymDefinitions_sorted.tex
\section*{Acronyms}

\acrodefplural{RES}[RES's]{Renewable Energy Sources}
\acrodefplural{AS}[AS's]{Ancillary Services}

\begin{acronym}[swissgrid]
\acro{ACE}{Area Control Error}
\acro{AGC}{Automatic Generation Control}
\acro{AHU}{Air Handling Unit}
\acro{AMI}{Advanced Metering Infrastructure}
\acro{AMR}{Advanced Metering Reading}
\acro{AS}{Ancillary Service}
\acrodefplural{AS}[AS]{Ancillary Services}
\acro{AVR}{Automatic Voltage Regulator}
\acro{BAS}{Building Automation System}
\acro{BESS}{Battery Energy Storage System}
\acro{BE}{Balancing Energy}
\acro{BG}{Balance Group}
\acro{CAISO}{California Independent System Operator}
\acro{CDF}{Cumulative Distribution Function}
\acro{COP}{Coefficient of Performance}
\acro{CWS}{Central Working Station}
\acro{DA}{Asymmetric Reserves at a Building with Daily Duration}
\acro{DAE}{Differential-Algebraic Equation}
\acro{DN}{Distribution Network}
\acro{DR}{Demand Response}
\acro{DSA}{Symmetric Reserves in Aggregate with Daily Duration}
\acro{DSB}{Symmetric Reserves at a Building with Daily Duration}
\acro{DSM}{Demand Side Management}
\acro{DSO}{Distribution System Operator}
\acro{EKZ}{the utility of the Kanton Zurich}
\acro{ENTSO-E}{European Network of Transmission System Operators for Electricity}
\acro{ERCOT}{Electric Reliability Council of Texas}
\acro{EWH}{Electric Water Heater}
\acro{FIT}{Feed-in Tariff}
\acro{FLEXLAB}{Facility for Low Energy eXperiments}
\acro{GPU}{Graphics Processing Unit}
\acro{HA}{Asymmetric Reserves at a Building with Hourly Duration}
\acro{HP}{Heat Pump}
\acro{HSA}{Symmetric Reserves in Aggregate with Hourly Duration}
\acro{HSB}{Symmetric Reserves at a Building with Hourly Duration}
\acro{HVAC}{Heating, Ventilation and Air-Conditioning}
\acro{HVDC}{High-Voltage Direct Current}
\acro{HV}{High-Voltage}
\acro{IRA}{Integrated Room Automation}
\acro{KiBaM}{Kinetic Battery Model}
\acro{LBNL}{Lawrence Berkeley National Laboratory}
\acro{LP}{Linear Program}
\acro{LV}{Low-Voltage}
\acro{MAE}{Mean Absolute Error}
\acro{MAPE}{Mean Absolute Percentage Error}
\acro{MHSE}{Moving Horizon State Estimation}
\acro{MILP}{Mixed-Integer Linear Program}
\acro{MLD}{Mixed Logical Dynamical}
\acro{MPC}{Model Predictive Control}
\acro{MPPT}{Maximum Power Point Tracking}
\acro{MV}{Medium-Voltage}
\acro{NPV}{Net Present Value}
\acro{ODE}{Ordinary Differential Equation}
\acro{OLTC}{On-Load Tap-Changer}
\acro{OPF}{Optimal Power Flow}
\acro{PC}{Power Constraints}
\acro{PCC}{Point of Common Coupling}
\acro{PDE}{Partial Differential Equation}
\acro{PDF}{Probability Distribution Function}
\acro{PEC}{Power and Energy Constraints}
\acro{PEV}{Plug-in Electric Vehicle}
\acro{PF}{Power Flow}
\acro{PFC}{Primary Frequency Control}
\acro{PI}{Proportional-Integral}
\acro{PID}{Proportional-Integral-Derivative}
\acro{PJM}{Pennsylvania, Jersey, and Maryland Power Market}
\acro{PLC}{Power Line Communication}
\acro{PSS}{Power System Stabilizer}
\acro{PV}{Photovoltaics}
\acro{QP}{Quadratic Program}
\acro{RBC}{Rule-based Control}
\acro{RC}{Resistance-capacitance}
\acro{REF}{Refrigerator}
\acro{RES}{Renewable Energy Source}
\acrodefplural{RES}[RES]{Renewable Energy Sources}
\acro{RHC}{Receding Horizon Control}
\acro{RMSE}{Root Mean Squared Error}
\acro{RPC}{Robust Problem with Power Constraints}
\acro{RPEC}{Robust Problem with Power and Energy Constraints}
\acro{SAT}{Supply Air Temperature}
\acro{SC}{Slab Cooling}
\acro{SFC}{Secondary Frequency Control}
\acro{SLP}{Sequential Linear Programming}
\acro{SM}{Smart Meter}
\acro{SOC}[SoC]{State of Charge}
\acro{SOH}[SoH]{State of Health}
\acro{SPC}{Stochastic Problem with Power Constraints}
\acro{SPEC}{Stochastic Problem with power and Energy Constraints}
\acro{TABS}{Thermally Activated Building Systems}
\acro{TCL}{Thermostatically Controlled Load}
\acro{TSO}{Transmission System Operator}
\acro{sLP}[SLP]{Stochastic Linear Program}
\acro{swissgrid}{Swiss Transmission System Operator}
\acro{UL}{Uncontrollable Load}
\acro{VAV}{Variable Air Volume}
\acro{VFD}{Variable Frequency Drive}
\end{acronym}

%% file: introduction_p1.tex
\section{Introduction}
\subsection{Motivation and Related Work}
Power system frequency reflects the balance between generation and demand of electric power. If generation exactly meets demand, the frequency is at its nominal value ($50$ Hz in Europe and $60$ Hz in North America). On the other hand, if generation becomes lower than demand, the frequency drops and vice versa. \acp{TSO} rely on frequency control reserves in the form of \acp{AS} to stabilize frequency after a sudden disturbance and recover it to its nominal value.

The integration of fluctuating \acp{RES} in the grid increases the need for frequency control reserves \cite{makarov_operational_2009}. Although these reserves are traditionally provided by power plants, additional reserve resources will be needed with large RES shares. Conceptually, loads can provide frequency control by reducing their consumption when frequency is low and increasing consumption when frequency is high \cite{Schweppe1980}.

\ac{HVAC} systems in commercial buildings are well suited for frequency control for three main reasons: (i) commercial HVAC systems make up a large percentage of the total electricity demand of a country (around $20\%$ in the US \cite{Hao2013Allerton,building_data}), (ii) commercial buildings often have a large thermal inertia, and (iii) many buildings (one-third of all buildings in the US \cite{Hao2013Allerton}) have a \ac{BAS} that facilitates control implementation. However, HVAC systems are typically complex with many control variables and cascaded control loops. Most of the early work on commercial buildings focused on the development of building thermal models \cite{Crabb1987,Kramer2012}, and on using the building's thermal mass for load shifting and peak shedding to minimize electricity cost and demand charges \cite{Braun1990,Braun1994,Henze1997}.

Some works investigated the potential of commercial buildings for \acp{AS} provision. A retail store and an office building participated in a pilot program for non-spinning reserves in the California Independent System Operator's AS market using global temperature adjustments in \cite{kiliccote2010open}. In \cite{kirby2008spinning} spinning reserve with a duration of $15$ minutes was successfully offered by curtailing the air conditioning load of a hotel. Ref. \cite{blum2014dynamic} used a detailed model of a \ac{VAV} HVAC system to simulate the provision of spinning reserve with setpoint adjustments in zone temperature, duct static pressure, \ac{SAT}, and chilled water temperature.

This paper concerns frequency regulation from commercial buildings, which is also known as secondary (or load) frequency control, automatic generation control, and frequency restoration reserve. Frequency regulation is activated via a signal sent from the \ac{TSO} typically every $2-4$ seconds with the goal of correcting frequency and tie-line power deviations \cite{Kundur}. There is a limited amount of theoretical, simulation-based or experimental work on frequency regulation with commercial buildings. In \cite{Kawachi2011} a heat pump was controlled to track a frequency regulation signal by changing the refrigerant's flow rate. Adjustments of the duct static pressure setpoint were used in a simulation study in \cite{zhao2013evaluation} for frequency regulation. Refs. \cite{Hao2013Allerton} and \cite{Hao2013ACC} investigated frequency regulation via fan power control and simulations showed that $15\%$ of the fan power can be offered as reserve, when the frequency band of the regulation signal is $f \in [1/(10~\textrm{min}), 1/(4~\textrm{sec})]$. The follow-up work \cite{LinIEEESGCom2013} included chiller control enlarging the frequency band to $1/(60~\textrm{min})$.

Since buildings are energy-constrained resources it is important to determine the reserve capacity reliably. \ac{MPC} was used in \cite{MaasoumyACC2014paper2} to quantify the flexibility of a commercial building and offer it to a utility. Ref. \cite{VrettosIFAC2014} presented a hierarchical control framework consisting of a reserve capacity scheduler, an MPC for HVAC system control, and a feedback controller to track the regulation signal. The framework of \cite{VrettosIFAC2014} was extended in \cite{VrettosTPS2016} to include energy-constrained regulation signals, and in \cite{VrettosTSTE2016} with a chance-constrained reserve scheduling formulation. Ref. \cite{Henze2014} proposed a simulation-based approach to estimate the reserve capacity neglecting the time-coupling across different scheduling intervals. The energy capacity of a commercial building was estimated in \cite{hughes2016identification} with a virtual battery model.

Estimating the building's baseline consumption without frequency regulation is challenging. Baseline estimation was performed on-line in \cite{LinTSG2014} using a low-pass filter. If MPC is used as in \cite{MaasoumyACC2014paper2,VrettosIFAC2014,VrettosTPS2016,VrettosTSTE2016} the baseline power is known ahead of time, which is advantageous because it facilitates the financial settlement.

Apart from simulation-based studies, experimental verification of frequency regulation by commercial buildings is necessary to build confidence for wide-spread implementation. Unfortunately, there have been only a few demonstrations and field tests so far. The feasibility of offering up- and down-regulation products with global temperature adjustments and ventilation power control was investigated in \cite{kiliccote_field_2012}. Ref.~\cite{MacDonald2014PJM} demonstrated that fans can provide frequency regulation with open-loop control of the frequency of the \ac{VFD} using industry-standard demand response communications. In \cite{LinTSG2014} an auditorium of a university building provided frequency regulation controlling the fan speed and air flow rate setpoints. The fan power was indirectly controlled via static duct pressure setpoint adjustments in \cite{MaasoumyACC2014paper2}.

Frequency regulation experiments with a variable speed heat pump were reported in \cite{kim2016experimental} in a lab-scale microgrid, using direct compressor control and adjustments of the supply water temperature setpoint. Ref.~\cite{Su2015Demonstration_p1} developed a \ac{PID} controller for frequency regulation with a chiller, and combined it with a high-pass filter of the regulation signal and a baseline estimator. The follow-up work \cite{Su2015Demonstration_p2} identified the \ac{BAS} delays and chiller ramp-rates and minimum cooling power limits as important issues for practical implementation. Finally, \cite{Backhaus2015} investigated experimentally the efficiency of fast demand response actions in commercial buildings.

\subsection{Contribution and Organization of this Paper}
To the best of our knowledge, this two-part paper presents the first experimental demonstration of frequency regulation from a commercial building that simultaneously addresses the following challenges: (i) a priori determination of reserve capacity and bidding in a day-ahead \ac{AS} market; (ii) a priori declaration of the short-term operating power around which we provide frequency regulation; (iii) balancing energy consumption and reserve capacity, such that the net profit is maximized and the effect on occupant comfort is minimal; and (iv) accurate tracking of the regulation signal with fan speed control.

To this end, we adopt the hierarchical control framework of \cite{VrettosIFAC2014,VrettosTPS2016,VrettosTSTE2016} and extend it in several ways. We determine the reserve capacity in a day-ahead fashion using robust optimization, and control the HVAC system with a robust \ac{MPC}. In contrast to \cite{VrettosIFAC2014,VrettosTPS2016,VrettosTSTE2016}, we investigate the use of different building models in the reserve scheduler and MPC controller, and develop a Kalman filter to provide state estimates to the controller. Whereas \cite{VrettosIFAC2014,VrettosTPS2016,VrettosTSTE2016} considered water-based HVAC systems, we focus on VAV HVAC systems with fans and propose computationally tractable reformulations to account for the fan nonlinear dynamics. Furthermore, in contrast to \cite{MacDonald2014PJM} that used open-loop control and to \cite{LinTSG2014} that developed a standard \ac{PI} controller for the fan, we propose a switched controller with a feedforward and a feedback loop in order to track the regulation signal. The proposed controller demonstrates a high performance without compromising stability.

In Part I of this two-part paper, we introduce the test facility in Section~\ref{sec_facility}, identify simple building models in Section~\ref{sec_modeling}, and present the hierarchical control design in Sections~\ref{CtrlDesign_lv1} - \ref{CtrlDesign_lv3}. Extensive experimental results are reported in Part II \cite{VrettosTSG2016Exp_p2}.

%% file: expFacility_v3.tex
\section{Test Facility} \label{sec_facility}

\subsection{FLEXLAB: our Test Facility}
\begin{figure}[t]
\centering \includegraphics[width=0.95\textwidth]{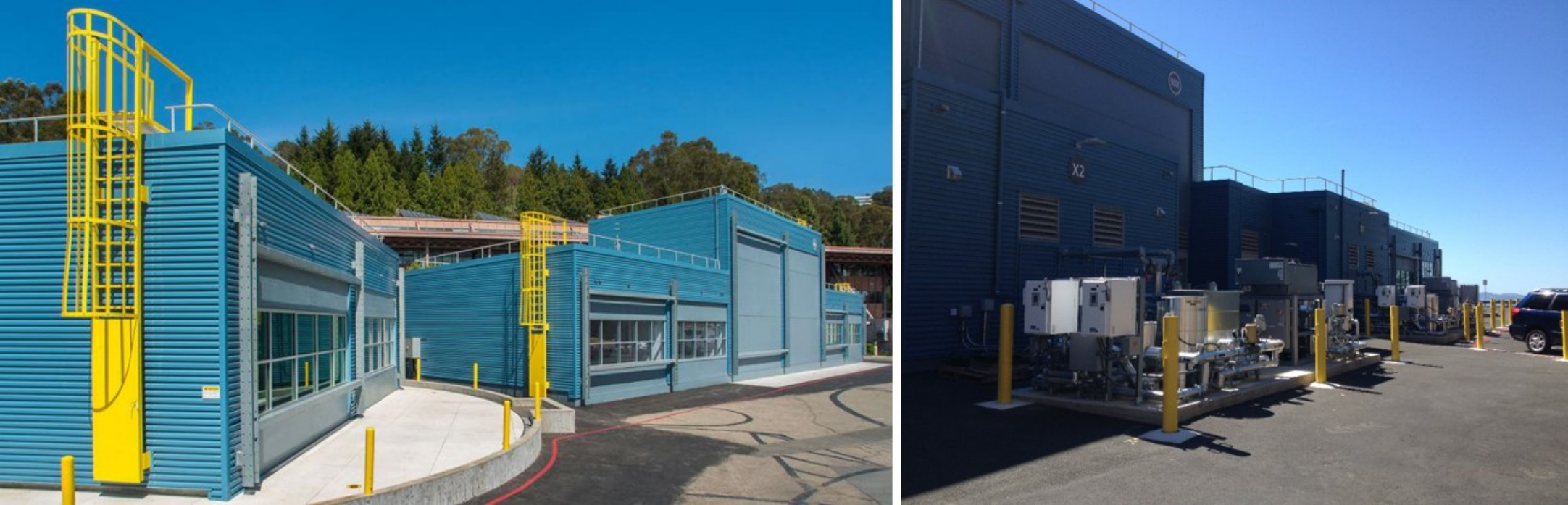}
\caption{The Facility for Low Energy eXperiments (\ac{FLEXLAB}) at LBNL.}
\label{fig:flexlab_pic}
\end{figure}

The experiment was performed at the Facility for Low Energy eXperiments (FLEX-LAB), a new facility for energy efficiency research in buildings located at the Lawrence Berkeley National Laboratory (LBNL). The facility (shown in Fig.~\ref{fig:flexlab_pic}) is comprised of $4$ buildings (called ``bays'') and each of them has $2$ thermally isolated test ``cells''. Each pair of cells is designed to be thermally identical, constructed with the same materials and dimensions.

The thermal isolation resulting from the near adiabatic walls between the two cells allows them to be modeled independently. This is a unique feature of FLEX-LAB that allows us to perform frequency regulation experiments in one of the two identical cells (cell ``1A''), while using the other one (cell ``1B'') as a benchmark to evaluate the effect of our control actions in real time and under the same external conditions. The bay used in our experiment has a south orientation and a total floor area of $120~\text{m}^2$ ($60~\text{m}^2$ per building cell).

\begin{figure}[t]
\centering \includegraphics[width=0.95\textwidth]{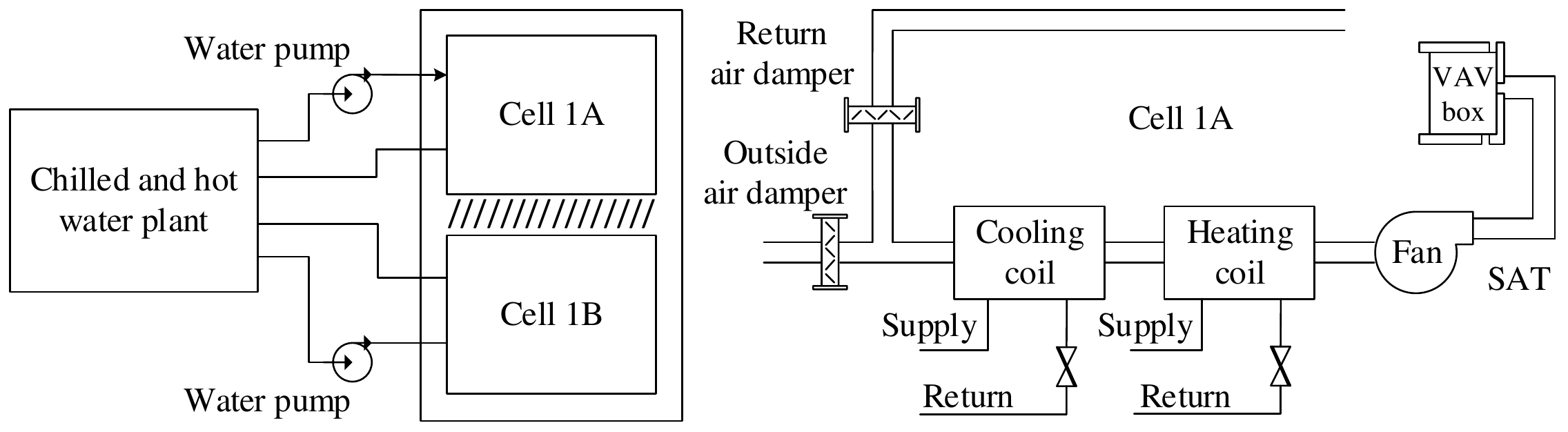}
\caption{The HVAC system of building cells 1A and 1B of FLEXLAB.}
\label{fig:HVAC_plot}
\end{figure}

Three major cascade control loops are present in a VAV HVAC system: chilled water temperature control, \ac{SAT} control, and zone (room) temperature control. A chiller plant cools down water that is then piped to the building's \ac{AHU}. The chilled water decreases the temperature of a mixture of return and outside air in the AHU using a heat exchanger, and the flow of the chilled water is controlled to maintain a constant \ac{SAT}. The cooled air is circulated to the building zones through the duct system using a fan. The temperature of each zone is maintained close to the desired setpoint by controlling the damper position of the VAV box.

Typically, an AHU provides several building zones with cooled air, whereas reheating is performed at the VAV boxes. As shown in Fig.~\ref{fig:HVAC_plot}, FLEXLAB differs from this typical mode of operation in two ways: (i) the cells are served by dedicated AHUs that contain a heating coil; and (ii) the air volume is controlled by fan speed alone rather than damper position. Although small, the test facility is a good representation of commercial buildings with \ac{VAV} systems constructed in the $1980$'s and is highly controllable.

\subsection{Control Approach} \label{sec_comm_arch}
FLEXLAB is controlled by a \ac{CWS} based on an existing control sequence programmed in LabVIEW. From a TestStand National Instruments user interface, the operator can monitor the system and modify the setpoints of various control loops. We develop the hierarchical controller for frequency regulation externally in order to minimize potential conflicts with the LabVIEW code, and send the control commands to the CWS via a scripting environment.

\begin{figure}[t]
\centering
\includegraphics[width=0.95\textwidth]{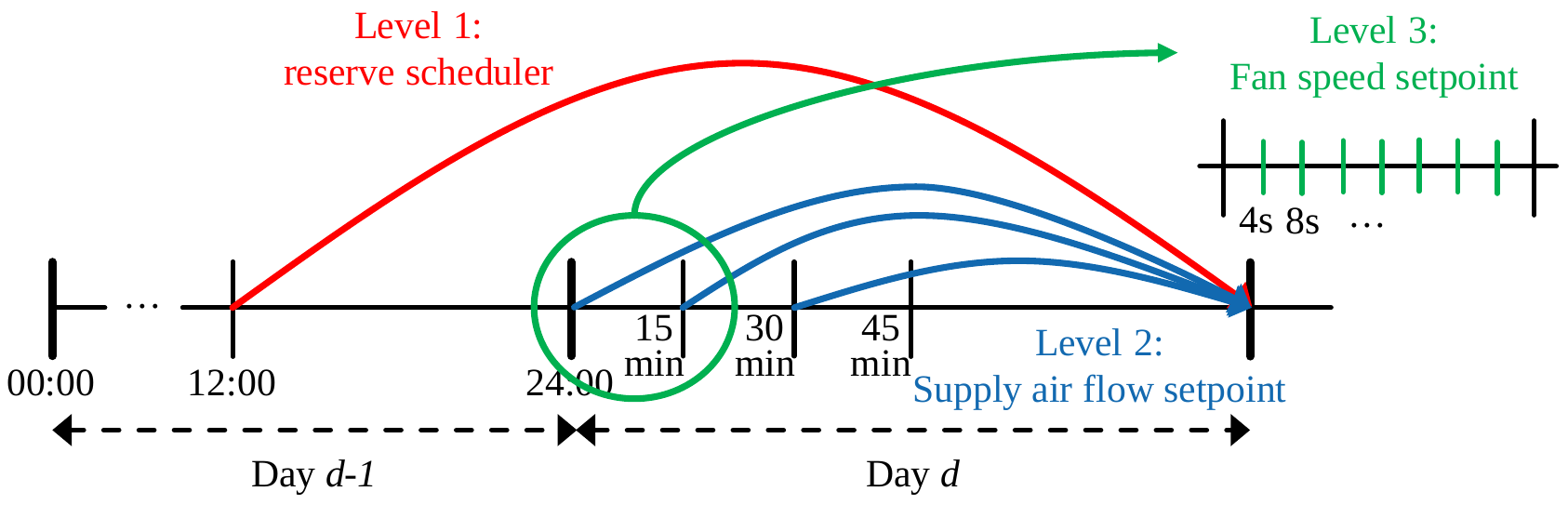}
\caption{Control sequences of the three levels of the hierarchical controller.} \label{fig:control_sequences}
\end{figure}

Specifically, we disable the zone temperature \ac{PI} control of FLEXLAB and replace it with an MPC-based controller, which determines the air flow rate setpoints. However, the hierarchical controller does not substitute the chilled water temperature and \ac{SAT} control loops of the HVAC system, which remain active. The hierarchical controller consists of the following three levels, and the control sequences are shown in Fig.~\ref{fig:control_sequences}.

\subsubsection{Level $1$: Reserve Scheduler}
The goal of the reserve scheduler is to determine the reserve capacity that the building can reliably offer to the \ac{TSO} by solving a multi-period robust optimization problem. We assume a day-ahead reserve scheduling occurring at $12.00$ of each day to determine the reserve capacity for the next day, which is common in several AS markets \cite{ENTSOE_AS_survey2014}.

\subsubsection{Level $2$: Room Climate Controller}
This zonal controller calculates the supply air flow rate setpoints that minimize energy consumption while ensuring occupant comfort under reserve provision. It is implemented as a robust MPC that runs every $15$ minutes along with a Kalman filter.

\subsubsection{Level $3$: Frequency Regulation Controller}
The goal of this controller is to track the frequency regulation signal every $4$ seconds by modifying the fan power via fan speed control with a \ac{VFD}. For this purpose, we designed a switched controller comprised of a feedforward model-based controller and a feedback \ac{PI} controller.

%% file: modeling_v3.tex
\section{Modeling and Identification} \label{sec_modeling}

\subsection{Building Thermal Model}
\begin{figure}[t]
\centering
\includegraphics[width=0.75\textwidth]{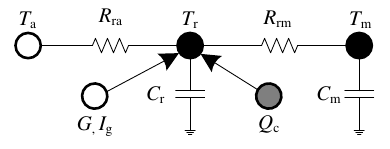}
\caption{The resistance-capacitance network of the building thermal model.} \label{fig:buildingRC}
\end{figure}

We model the building with the 2-state resistance-capacitance network of Fig.~\ref{fig:buildingRC}. If the heating coil of the AHU is deactivated, the cooling power of the HVAC system is given by $Q_\text{c}=\dot{m} c_\text{p} (T_\text{s} - T_\text{r})$, where $T_\text{r}$ is the room temperature, $\dot{m}$ is the mass air flow rate, $c_\text{p}$ is the specific heat capacity of air, and $T_\text{s}$ is the \ac{SAT}. Let us denote by: $T_\text{m}$ the temperature of the building's lumped thermal mass; $C_\text{r}$ and $C_\text{m}$ the thermal capacitances of the room and the thermal mass; $T_\text{a}$ the ambient temperature; $R_\text{ra}$ and $R_\text{rm}$ the thermal resistances between the room and the ambient, and between the room and the thermal mass, respectively; $G$ the solar irradiance; $\gamma$ the solar irradiance absorption factor; and $I_\text{g}$ the internal heat gain. With this notation, the continuous-time state-space model can be written as
\begin{align} \label{ssCon}
\begin{bmatrix}
\dot{T}_\text{r} \\ \dot{T}_\text{m}
\end{bmatrix} &=
\overbrace{\begin{bmatrix*}[l]
-\left(\frac{1}{C_\text{r} R_\text{ra}}+\frac{1}{C_\text{r} R_{\text{rm}}}\right) \quad \frac{1}{C_\text{r} R_{\text{rm}}} \\ \frac{1}{C_\text{m} R_{\text{rm}}} \quad\quad\quad\quad\quad-\frac{1}{C_\text{m} R_{\text{rm}}}
\end{bmatrix*}}^{A^c} \cdot
\begin{bmatrix}
T_\text{r} \\ T_\text{m}
\end{bmatrix} +
\overbrace{\begin{bmatrix}
\frac{c_\text{p} T_\text{s}}{C_\text{r}} \\ 0
\end{bmatrix}}^{B_u^c} \dot{m} + \nonumber \\
&\overbrace{\begin{bmatrix*}[l]
-\frac{c_\text{p}}{C_\text{r}} ~~ 0 \\ 0 \quad\quad 0
\end{bmatrix*}}^{B_{xu}^c} \cdot \begin{bmatrix} T_\text{r} \\ T_\text{m} \end{bmatrix} \dot{m} +
\overbrace{\begin{bmatrix*}[l]
\frac{1}{C_\text{r} R_\text{ra}} ~~ \frac{\gamma}{C_\text{r}} ~~ \frac{1}{C_\text{r}} \\
0 ~~\quad\quad 0 ~~\quad 0
\end{bmatrix*}}^{B_v^c} \cdot
\begin{bmatrix}
T_\text{a} \\ G \\ I_\text{g}
\end{bmatrix}.
\end{align}
The model is bilinear between the control input $\dot{m}$ and the state $T_\text{r}$. Note that there is no bilinearity between $\dot{m}$ and $T_\text{s}$ because $T_\text{s}$ is fixed in our experiment.

With a first-order Euler discretization, the discrete-time model maintains the structure of the continuous-time matrices \cite{Siroky2011}
\begin{align} \label{ssDis}
&\begin{bmatrix}
T_{\text{r},k+1} \\ T_{\text{m},k+1}
\end{bmatrix} =
\overbrace{
\begin{bmatrix*}[l]
a_{11} ~~ a_{12} \\ a_{21} ~~ a_{22}
\end{bmatrix*}
}^{A} \cdot
\begin{bmatrix}
T_{\text{r},k} \\ T_{\text{m},k}
\end{bmatrix} +
\overbrace{
\begin{bmatrix}
b T_{\text{s}} \\ 0
\end{bmatrix}
}^{B_{u}}
\dot{m}_k + \nonumber \\
&\quad\quad\overbrace{
\begin{bmatrix*}[l]
-b ~~ 0 \\ 0 \quad 0
\end{bmatrix*}
}^{B_{xu}} \cdot
\begin{bmatrix} T_{\text{r},k} \\ T_{\text{m},k} \end{bmatrix} \dot{m}_k +
\overbrace{
\begin{bmatrix*}[l]
d_{11} ~~ d_{12} ~~ d_{13} \\
0 \quad\quad 0 \quad\quad 0
\end{bmatrix*}
}^{B_{v}} \cdot
\begin{bmatrix}
T_{\text{a},k} \\ G_k \\ I_{\text{g},k}
\end{bmatrix},
\end{align}
where the state, input and disturbance vectors are defined as
\begin{align}
x_k = [T_{\text{r},k} ~~ T_{\text{m},k}]^\top,~~
u_k = \dot{m}_k,~~
v_k = [T_{\text{a},k} ~~ G_k ~~ I_{\text{g},k}]^\top~.
\end{align}

We developed the following regression problem to identify the entries of matrices $A$, $B_u$, $B_{xu}$, and $B_v$ using measurements of $T_{\text{r},k}$, $\dot{m}_k$, $T_{\text{a},k}$, $G_k$, and $I_{\text{g},k}$
\begin{subequations} \label{modelIDproblem}
\begin{align}
    &\underset{A,B_u,B_{xu},B_v,\hat{x}_k}{\min} \sum_k \left[\hat{x}_k(1) - x_k(1)\right]^2 \\
    &\text{s.t.~} \hat{x}_{k+1} = A x_k + B_{u} u_k + B_{xu} x_k u_k + B_v v_k,~\forall k\label{dynCons}\\
    &\text{~~~~}    a_{12},~a_{21},~b,~d_{11},~d_{12},~d_{13} \geq 0 \label{structureCons} \\
    &\text{~~~~}    \hat{T}^\text{min}_{\text{m},k} \leq \hat{x}_k(2)  \leq \hat{T}^\text{max}_{\text{m},k} \label{T2Cons},~\forall k \\
    &\text{~~~~}    |\text{eig}(A)| \leq 1,~ |\text{eig}\big(A+B_{xu} u_k\big)| \leq 1,~\forall k \label{stabCons}~,
\end{align}
\end{subequations}
where $x_k(1)=T_{\text{r},k}$ and $x_k(2)=T_{\text{m},k}$. Since the state $T_{\text{m},k}$ is not directly measured, the regression \eqref{modelIDproblem} is a non-linear optimization problem that involves multiplications of the optimization variables (the model parameters).

Constraints \eqref{structureCons} represent the fact that the positive elements of the continuous-time matrices $A^c$, $B_u^c$, and $B_v^c$ in \eqref{ssCon} remain positive in the the discrete-time matrices $A$, $B_u$, and $B_v$ due to the first-order discretization. Indeed, the discrete-time matrices are computed with $A = I+ \Delta{t} \cdot A^c$, $B_u=\Delta{t} \cdot B_u^c$, and $B_v=\Delta{t} \cdot B_v^c$, where $\Delta{t}=15$~minutes is the discretization step and $I$ is the identity matrix. Constraints \eqref{T2Cons} are lower and upper bounds on the estimated unmeasured state $\hat{x}_k(2)=\hat{T}_{\text{m},k}$ (the bounds $\hat{T}^\text{min}_{\text{m},k}=0.01 \cdot T_{\text{r},k}$ and $\hat{T}^\text{max}_{\text{m},k}=2.5 \cdot T_{\text{r},k}$ were used). Constraints \eqref{stabCons} impose stability of the identified bilinear model to avoid over-fitting if the data set is small or if it does not have sufficient excitation in terms of control inputs.

\begin{figure}[t]
\centering
\includegraphics[width=0.99\textwidth]{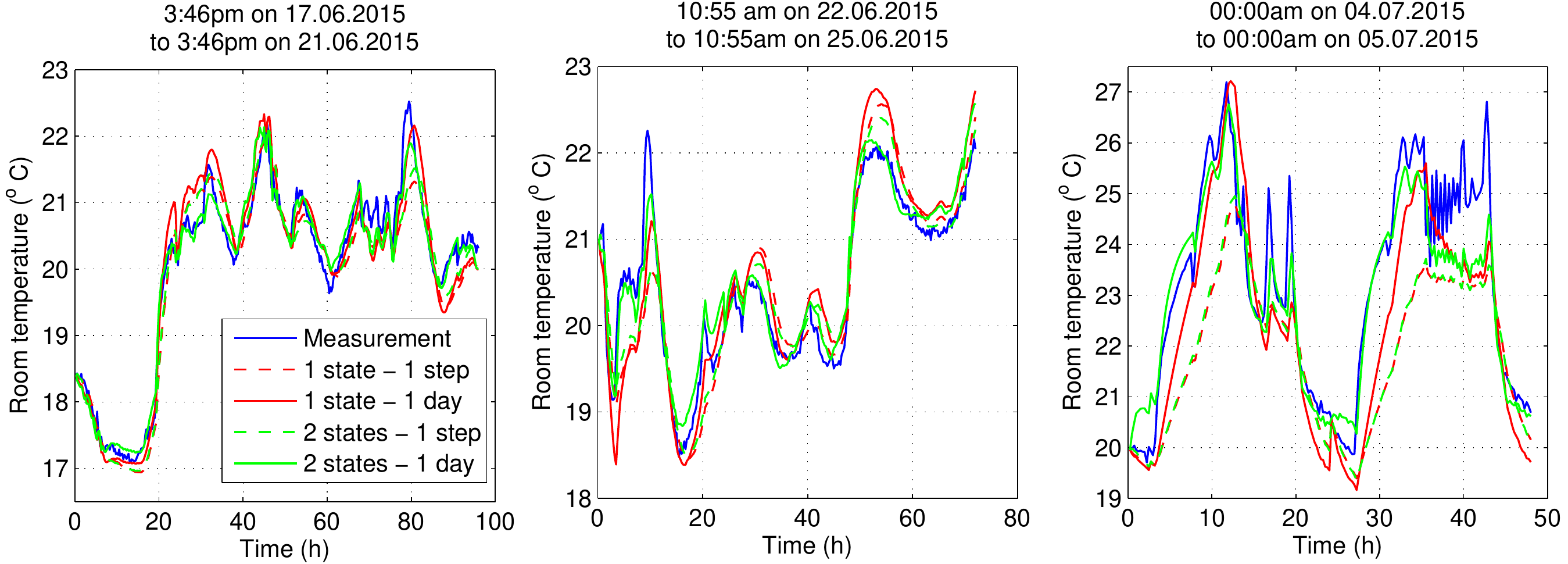}
\caption{Identification and comparison of four building thermal models.} \label{fig:thermal_id_results}
\end{figure}

\begin{table}[h!]
\renewcommand{\arraystretch}{1.1}
\caption{Comparison of building models with respect to RMSE}
\centering
\begin{tabular}{c|c|c|c|c}
\hline
& $1$~state, $1$~step & $1$~state, $1$~day & $2$~states, $1$~step & $2$~states, $1$~day \\
RMSE & $0.92^\circ$C & $0.67^\circ$C & $0.89^\circ$C & $0.42^\circ$C\\
\hline
\end{tabular}
\label{tab:thermal_model_cmp}
\end{table}

\begin{table}[t!]
\renewcommand{\arraystretch}{1.1}
\caption{Parameters of the older model ($2$ states, $1$ day prediction)}
\centering
\begin{tabular}{cccc}
\hline
$a_{11}=0.8665$ & $a_{12}=0.0918$ & $a_{21}=0.0374$ & $a_{22}=0.9703$ \\
$b=0.2996$ & $d_{11}=0.0230$ & $d_{12}=2.016 \cdot 10^{-4}$ & $d_{13}=1.424 \cdot10^{-4}$ \\
\hline
\end{tabular}
\label{tab:old_thermal_model_params}
\end{table}

\begin{table}[t!]
\renewcommand{\arraystretch}{1.1}
\caption{Parameters of the new model ($2$ states, $1$ day prediction)}
\centering
\begin{tabular}{cccc}
\hline
$a_{11}=0.6344$ & $a_{12}=0.2661$ & $a_{21}=0.1021$ & $a_{22}=0.9170$ \\
$b=0.4716$ & $d_{11}=0.0405$ & $d_{12}=0.0028$ & $d_{13}=3.3686 \cdot10^{-4}$ \\
\hline
\end{tabular}
\label{tab:recent_thermal_model_params}
\end{table}

\begin{table}[t!]
\renewcommand{\arraystretch}{1.1}
\caption{Fan model parameters}
\centering
\begin{tabular}{cccc}
\hline
$\alpha_3=2588.2$ & $\alpha_2=-1458.0$ & $\alpha_1=630.9$ & $\alpha_0=28.7$\\
$\beta_3=0.0032$ & $\beta_2=-0.0151$ & $\beta_1=1.4521$ & $\beta_0=55.7634$ \\
 -  & - & $\gamma_1=0.0133$ & $\gamma_0=0.0606$ \\
\hline
\end{tabular}
\label{tab:fan_model_params}
\end{table}

The optimization problem \eqref{modelIDproblem} identifies a model based on ``$1$-step ahead prediction'' meaning that the estimate $\hat{x}_{k+1}(1)=\hat{T}_{\text{r},k+1}$ at time step $k+1$ depends directly on the measurement $x_k(1)=T_{\text{r},k}$ at time step $k$. This model identification approach is standard for MPC applications because the one-step prediction is important. Nevertheless, in this experiment the building model is also used in the reserve scheduler, and thus high-quality day-ahead predictions are also important.

For this purpose, we propose to modify the regression problem by substituting the measurement $x_k(1)=T_{\text{r},k}$ with the optimization variable $\hat{x}_k(1)=\hat{T}_{\text{r},k}$ in \eqref{dynCons}. We use this more complex formulation to obtain a ``$1$-day ahead prediction'' model, which is expected to predict the building states up to one day ahead more accurately, because it is more flexible than the original formulation. Of course, it is also possible to identify a $1^\text{st}$-order model by neglecting the lumped thermal mass of the room, which simplifies the regression problem significantly.

Two sets of building model parameters were fitted to investigate the importance of periodic calibration. The first set (``older model'') used data from $17-25$ June and $4-5$ July $2015$, whereas the second set (``new model'') used data from $12-18$ November $2015$. The building was excited with different combinations of air flow rate, SAT, and internal heat gains. Four different model variants were compared: (i) $1$-state model with $1$-step ahead prediction, (ii) $1$-state model with $1$-day ahead prediction, (iii) $2$-state model with $1$-step ahead prediction, and (iv) $2$-state model with $1$-day ahead prediction.

The identification results are shown in Fig.~\ref{fig:thermal_id_results} and the model \acp{RMSE} are given in Table~\ref{tab:thermal_model_cmp}. As expected, increasing the number of states or the prediction horizon reduces the RMSE. We use the $2$-state model with $1$-day ahead prediction in the frequency regulation experiments because it achieves the lowest RMSE. The identified model parameters are shown in Tables~\ref{tab:old_thermal_model_params} and \ref{tab:recent_thermal_model_params}.

\subsection{Fan Model} \label{fan_model}
A steady-state fan model is required in the MPC to map the air flow setpoint to fan power, and in the frequency regulation controller to convert the electric power setpoint to a fan speed reference. According to the fan laws, the mass air flow rate $u$ is proportional to the fan speed $N_\text{f}$, and the fan power $P_\text{f}$ increases with the cube of the fan speed. Therefore, a steady-state model can be obtained by fitting the parameters of
\begin{align}
P_\text{f} &= f(u) = \alpha_3 u^3 + \alpha_2 u^2 + \alpha_1 u + \alpha_0 \label{flow-to-power}\\
P_\text{f} &= g(N_\text{f}) = \beta_3 N_\text{f}^3 + \beta_2 N_\text{f}^2 + \beta_1 N_\text{f} + \beta_0 \label{speed-to-power} \\
u &= h(N_\text{f}) = \gamma_1 N_\text{f} + \gamma_0~. \label{speed-to-flow}
\end{align}

\begin{figure}[t]
\centering
\includegraphics[width=0.95\textwidth]{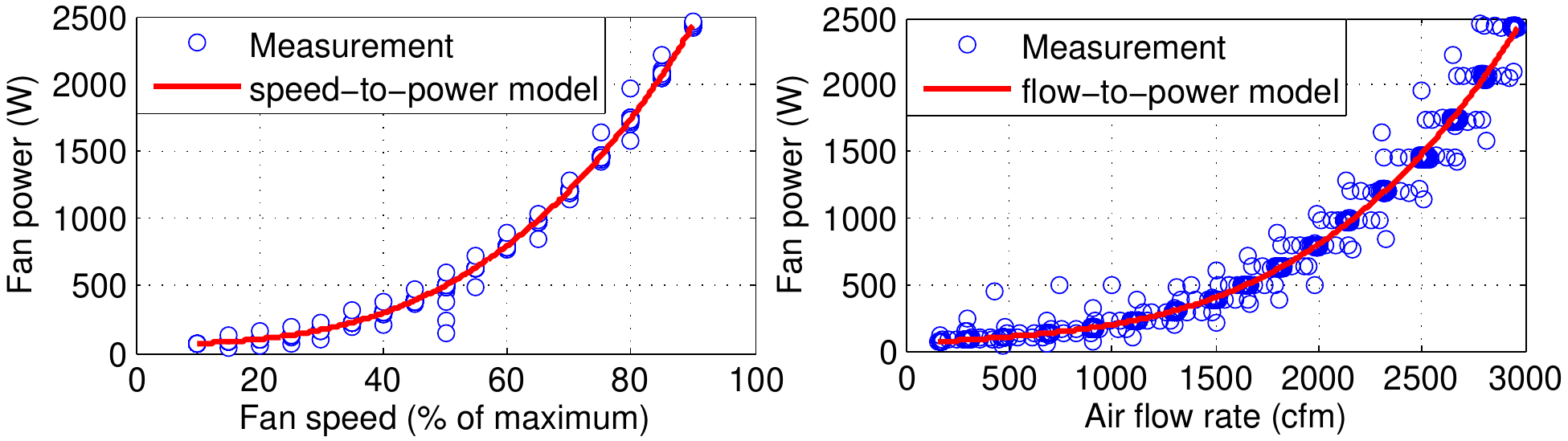}
\caption{Raw fan measurements and the identified fan models.} \label{fig:fan_id_results}
\end{figure}

For this purpose, we vary the fan speed setpoint from the minimum value of $10\%$ to the maximum value of $90\%$ of the rated fan power (with a step of $5\%$) and record the air flow rate and electric power. Each setpoint is kept for $6$ minutes, but the first $20$ seconds of the data after each step change are discarded to account for communication delays and the fan transients. The identified parameters are given in Table~\ref{tab:fan_model_params}, whereas the measurements and identified models are shown in Fig.~\ref{fig:fan_id_results}. The fitting performance is very high: the \ac{RMSE} is only $5$ W for the speed-to-power model and $21$ W for the flow-to-power model.

%% file: controlDesign_v3.tex
\section{Level $1$: Reserve Scheduler} \label{CtrlDesign_lv1}

\subsection{Robust Reserve Scheduling Formulation}
\begin{figure}[t]
\centering \includegraphics[width=0.85\textwidth]{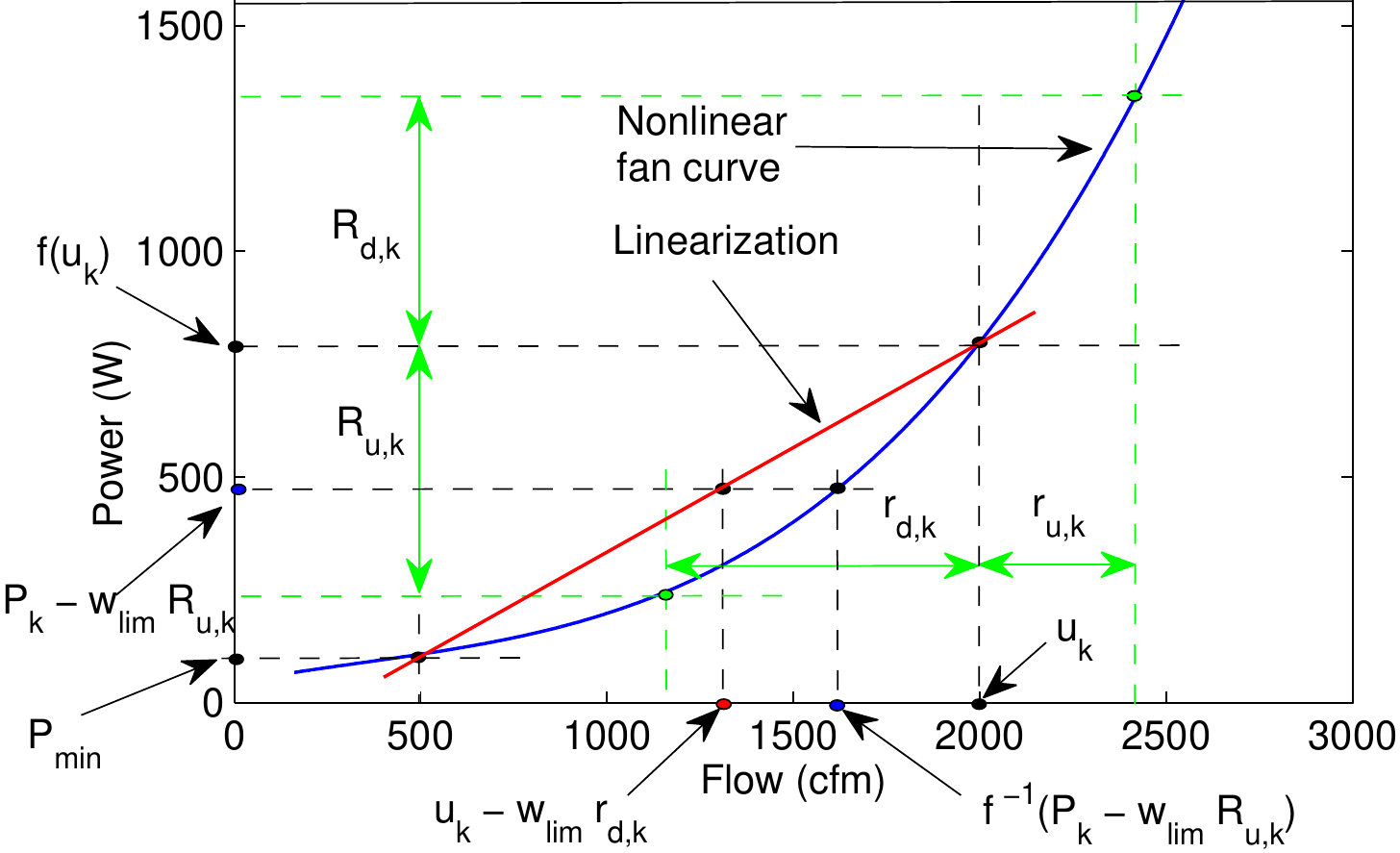}
\caption{The fan curve linearization for the optimization. The thermal reserves and the electric reserves are shown around an air flow rate operating point.}
\label{fig:Fan_model_fit_explain}
\end{figure}

Let $R_{\text{u},k}$ and $R_{\text{d},k}$ denote the electric reserve capacities at time step $k$ for regulation up and down, respectively. \footnote{Up-reserve is an increase of generation or decrease of consumption, whereas down-reserve is a decrease of generation or increase of consumption.} It is convenient to define also the thermal up- and down-reserve capacities $r_{\text{u},k}$ and $r_{\text{d},k}$ as the maximum changes in the mass air flow rate due to reserve provision. In cooling operation, a request for regulation up results in a reduction in air mass flow rate, such that $R_{\text{u},k}$ is related to $r_{\text{d},k}$. On the other hand, regulation down results in an increase in air mass flow rate ($R_{\text{d},k}$ is related to $r_{\text{u},k}$). $R_{\text{u},k}$ and $R_{\text{d},k}$ are coupled to $r_{\text{d},k}$ and $r_{\text{u},k}$ with the flow-to-power fan model \eqref{flow-to-power} according to
\begin{align}
R_{\text{u},k} &= f(u_k)-f(u_k-r_{\text{d},k}) \label{Ru_def}\\
R_{\text{d},k} & =f(u_k+r_{\text{u},k})-f(u_k) \label{Rd_def}~,
\end{align}
where $u_k$ is the operating point of air flow. This nonlinear relationship is graphically shown in Fig~\ref{fig:Fan_model_fit_explain}.

The objective is to minimize the total cost defined as the sum of energy consumption cost and reserve profit
\begin{align}
c_k P_k + \lambda_k \left(R_{\text{d},k} + R_{\text{u},k}\right)~,
\end{align}
where $c_k$ is the electricity price and $P_k$ is the fan power consumption. Assuming the same payment $\lambda_k$ for up- and down-reserves and using \eqref{Ru_def} and \eqref{Rd_def}, the reserve profit is given by
\begin{align}
\lambda_k \left(R_{\text{d},k} + R_{\text{u},k}\right) &= \lambda_k \big[f(u_k+r_{\text{u},k})-f(u_k-r_{\text{d},k})\big].
\end{align}

Typically, the \ac{TSO} requests the reserve energy as a fraction of the reserve capacity using a normalized frequency regulation signal $w_k \in [-1,1]$ \cite{swissgridAS}. Thus, the reserve at time step $k$ is
\begin{align} \label{res_energy_def}
R_k = \begin{cases}
        w_k R_{\text{u},k}, & \textrm{if} \hspace{2mm} w_k<0\\
        w_k R_{\text{d},k}, & \textrm{if} \hspace{2mm} w_k \geq 0~.
    \end{cases}
\end{align}
The electric reserve request can be translated to a perturbation around $u_k$ using the fan curve
\begin{align} \label{Du_def}
\Delta{u}_k = f^{-1}\left(P_k + R_k \right) - u_k~.
\end{align}

With the above notation, the multi-period robust reserve scheduling problem can be written as
\begin{subequations} \label{Lv1}
\begin{align}
	&\min_{u_k,r_{\text{u},k},r_{\text{d},k}} \sum\nolimits_{k=0}^{N_1-1} c_k f(u_k) -\lambda_k \left(R_{\text{d},k} + R_{\text{u},k}\right) \label{OFLv1}\\
    &\text{s.t.~}  x_{k+1} = A x_k + B_u T_{\text{s}} \cdot (u_k+\Delta{u}_k) + \nonumber\\
     &\quad~~\quad~~\quad~~ B_{xu} x_k \cdot (u_k+\Delta{u}_k) + B_v v_k,~\forall k \label{DynamicsLv1} \\
     &\quad~~u_{\text{min},1} \leq u_k+\Delta{u}_k \leq u_{\text{max},1}, ~\forall w_k \in [-1,1],~\forall k \label{InputConLv1} \\
     &\quad~~x_{\text{min},k} \leq x_k \leq x_{\text{max},k}, ~\forall w_k \in [-w_\text{lim},w_\text{lim}],~\forall k~. \label{StateConLv1}
\end{align}
\end{subequations}
Equation \eqref{DynamicsLv1} represents the building dynamics, whereas \eqref{InputConLv1} and \eqref{StateConLv1} set upper and lower bounds on the air mass flow rate and temperature, respectively. The limits $u_{\text{min},1}$ and $u_{\text{max},1}$ are calculated at the fan speeds $20\%$ and $80\%$ with
\begin{align} \label{mdot_limits_lv1}
u_{\text{min},1} = h\left(20\%\right), \quad u_{\text{max},1} = h\left(80\%\right)~.
\end{align}
The comfort zone (temperature) limits $x_{\text{min},k}$ and $x_{\text{max},k}$ are time-varying and different for working and non-working hours.

Uncertainty in the unknown frequency regulation signal, $w_k$, is explicitly handled in the reserve scheduling problem. Formulation \eqref{Lv1} builds robustness to $w_k$ with the robust input and state constraints \eqref{InputConLv1} and \eqref{StateConLv1}. Since the worst case in terms of power is either full up-reserve or down-reserve activation, $w_k$ can take any value in $[-1,1]$ in \eqref{InputConLv1}. The energy content of the regulation signal over $15$ minutes is typically limited, which is captured by the uncertainty set $w_k \in [-w_\text{lim},w_\text{lim}]$ in \eqref{StateConLv1}, where $0<w_\text{lim}\leq1$ is the energy limit.

\subsection{Reformulation and Approximation}
Due to the uncertain variable $w_k$, problem \eqref{Lv1} is not directly solvable. However, we derive the robust counterpart problem \eqref{Lv1rc} by formulating the input and state constraints of \eqref{Lv1} only for the boundaries of the uncertainty $w_k$. In \eqref{Lv1rc}, $\overline{x}_k$ and $\underline{x}_k$ are the worst case higher and lower state trajectories, respectively. We show in Proposition~\ref{reformulation_proposition} that problems \eqref{Lv1rc} and \eqref{Lv1} are equivalent, if the building operates in cooling mode.

\begin{subequations} \label{Lv1rc}
\begin{align}
    &\min_{u_k,r_{\text{u},k},r_{\text{d},k}} \sum\nolimits_{k=0}^{N_1-1} c_k f(u_k) -\lambda_k \left(R_{\text{d},k} + R_{\text{u},k}\right) \label{OFLv1rc}\\
    &\text{s.t.~} \overline{x}_{k+1} = A \overline{x}_k + B_u T_{\text{s}} \cdot f^{-1}\left(P_k-w_\text{lim} R_{\text{u},k}\right) +\nonumber \\
    &\quad~~\quad~~~~ B_{xu} \overline{x}_k \cdot f^{-1}\left(P_k-w_\text{lim} R_{\text{u},k}\right) + B_v v_k,~\forall k \label{DynamicsLv1rcUp} \\
    &\quad~~ \underline{x}_{k+1} = A \underline{x}_k + B_u T_{\text{s}} \cdot f^{-1}\left(P_k+w_\text{lim} R_{\text{d},k}\right) + \nonumber \\ &\quad~~\quad~~~~B_{xu} \underline{x}_k \cdot f^{-1}\left(P_k+w_\text{lim} R_{\text{d},k}\right) + B_v v_k,~\forall k \label{DynamicsLv1rcDown} \\
     &\quad~~ u_{\text{min},1} \leq u_k-r_{\text{d},k},~~u_k+r_{\text{u},k} \leq u_{\text{max},1}, \forall k \label{InputConLv1rc} \\
     &\quad~~ x_{\text{min},k} \leq \underline{x}_k,\quad\quad\quad~~\overline{x}_k \leq x_{\text{max},k} ~\forall k~. \label{StateConLv1rc}
\end{align}
\end{subequations}

\begin{lemma}\label{convex_monotonic}
Function $f(u)$ is both monotonic and convex.
\end{lemma}

\begin{proof}
It is sufficient to show that the $1^\text{st}$ and $2^\text{nd}$ order derivatives of $f$ are non-negative for the parameters of Table~\ref{tab:fan_model_params}.
\end{proof}

\begin{lemma}\label{min_max_eq}
If $w_k \in [-w_\text{lim}, w_\text{lim}]$ with $0< w_\text{lim}\leq1$, the following statements are true:
\begin{align}
\underset{w_k}{\min} \left(u_k+\Delta{u}_k \right) &= f^{-1}\left(P_k-w_\text{lim} R_{\text{u},k}\right)~\text{for}~w_\text{lim}\leq1 \label{lemma2_p1}\\
\underset{w_k}{\max} \left(u_k+\Delta{u}_k \right) &= f^{-1}\left(P_k+w_\text{lim} R_{\text{d},k}\right)~\text{for}~w_\text{lim}\leq1 \label{lemma2_p2}\\
\underset{w_k}{\min} \left(u_k+\Delta{u}_k \right) &= u_k-r_{\text{d},k,}~\text{for}~w_\text{lim}=1 \label{lemma2_p3}\\
\underset{w_k}{\max} \left(u_k+\Delta{u}_k \right) &= u_k+r_{\text{u},k,}~\text{for}~w_\text{lim}=1~.\label{lemma2_p4}
\end{align}
\end{lemma}

\begin{proof}
The proof is given in Appendix~\ref{proof_lemma2}.
\end{proof}

\begin{assumption} \label{SAT_assumption}
We assume that the building operates in cooling mode by deactivating the heating coil of the AHU, and the SAT is controlled to a setpoint $T_{\text{s}}$ that satisfies
\begin{align}
T_{\text{s}} \leq x_{\text{min},k} \leq T_{\text{r},k},~\forall k~.
\end{align}
Therefore, increasing the air flow rate will always decrease the room temperature.
\end{assumption}

\begin{proposition}\label{reformulation_proposition}
Under Assumption~\ref{SAT_assumption}, optimization problems \eqref{Lv1} and \eqref{Lv1rc} are equivalent.
\end{proposition}

\begin{proof}
The proof is given in Appendix~\ref{proof_propo1}.
\end{proof}

The dynamics in \eqref{Lv1rc} involve the inverse of a polynomial combination of optimization variables and are complex. This is in contrast to the formulations of \cite{VrettosIFAC2014,VrettosTPS2016,VrettosTSTE2016} where the nonlinear fan dynamics were not considered. We propose to approximate \eqref{Lv1rc} by the simple linearization of the inverse function shown in Fig.~\ref{fig:Fan_model_fit_explain}, which leads to problem \eqref{Lv1rca}. As shown by Propositions~\ref{bound_proposition} and \ref{exactness_proposition}, problems \eqref{Lv1rc} and \eqref{Lv1rca} are equivalent only for the special case $w_\text{lim}=1$, but not in general ($0<w_\text{lim}\leq1$).

\begin{subequations} \label{Lv1rca}
\begin{align}
	&\min_{u_k,r_{\text{u},k},r_{\text{d},k}} \sum\nolimits_{k=0}^{N_1-1} c_k f(u_k) -\lambda_k \left(R_{\text{d},k} + R_{\text{u},k}\right) \label{OFLv1rca}\\
    &\text{s.t.~} \overline{x}_{k+1} = A \overline{x}_k + B_u T_{\text{s}} \cdot (u_k-w_\text{lim} \cdot r_{\text{d},k}) +\nonumber \\
    &\quad~~\quad~~\quad~~B_{xu} \overline{x}_k \cdot (u_k-w_\text{lim} \cdot r_{\text{d},k}) + B_v v_k,~\forall k \label{DynamicsLv1rcaUp} \\
    &\quad~~ \underline{x}_{k+1} = A \underline{x}_k + B_u T_{\text{s}} \cdot (u_k+ w_\text{lim} \cdot r_{\text{u},k}) + \nonumber \\ &\quad~~\quad~~\quad~~B_{xu} \underline{x}_k \cdot (u_k+ w_\text{lim} \cdot r_{\text{u},k}) + B_v v_k,~\forall k \label{DynamicsLv1rcaDown} \\
     &\quad~~ \eqref{InputConLv1rc}, \eqref{StateConLv1rc}~.\nonumber
\end{align}
\end{subequations}

\begin{proposition}\label{bound_proposition}
Let $\overline{x}_k^\star$ and $\underline{x}_k^\star$ denote the maximum and minimum state trajectories of the original problem \eqref{Lv1rc}. Furthermore, let $\overline{x}_k^*$ and $\underline{x}_k^*$ denote the maximum and minimum state trajectories obtained by \eqref{Lv1rca}. With $0<w_\text{lim}\leq1$, $\overline{x}_k^* \geq \overline{x}_k^\star$ and $\underline{x}_k^* \geq \underline{x}_k^\star$ hold for any time step $k$, i.e., the approximation \eqref{Lv1rca} overestimates the maximum and minimum room temperatures compared with the original problem \eqref{Lv1rc}.
\end{proposition}

\begin{proof}
The proof is given in Appendix~\ref{proof_propo2}.
\end{proof}

\begin{proposition}\label{exactness_proposition}
Problems \eqref{Lv1rc}, \eqref{Lv1rca} are equivalent if $w_\text{lim}=1$.
\end{proposition}

\begin{proof}
The proof follows directly by rewriting \eqref{DynamicsLv1rcUp} and \eqref{DynamicsLv1rcDown} using \eqref{lemma2_p3}, \eqref{lemma2_p4} and \eqref{aux_eq2_lemma2} from Lemma~\ref{min_max_eq}.
\end{proof}

Note that the overestimation of the maximum and minimum room temperature builds additional robustness to temperature excursions above $x_{\text{max},k}$, but it reduces robustness to temperature excursions below $x_{\text{min},k}$. This is desirable because the state trajectory will generally remain closer to $x_{\text{max},k}$ than $x_{\text{min},k}$ due to minimization of energy consumption cost in \eqref{OFLv1rca}.

Problem \eqref{Lv1rca} is a deterministic nonlinear optimization problem with cubic objective function, bilinear equality constraints and linear inequality constraints. Although this is a non-convex problem, it is possible to solve it in due time with the solver IPOPT due to its relatively small size (the problem was compiled with YALMIP \cite{lofberg2004yalmip}). The outcome of \eqref{Lv1rca} is the up-reserve $R_{\text{u},k}$ and down-reserve $R_{\text{d},k}$ capacity for each time slot of the scheduling horizon $k \in [0,N_1-1]$.

\subsection{Energy Limit of Regulation Signal} \label{selection_wlim}
\begin{figure}[t]
\centering
\includegraphics[width=0.95\textwidth]{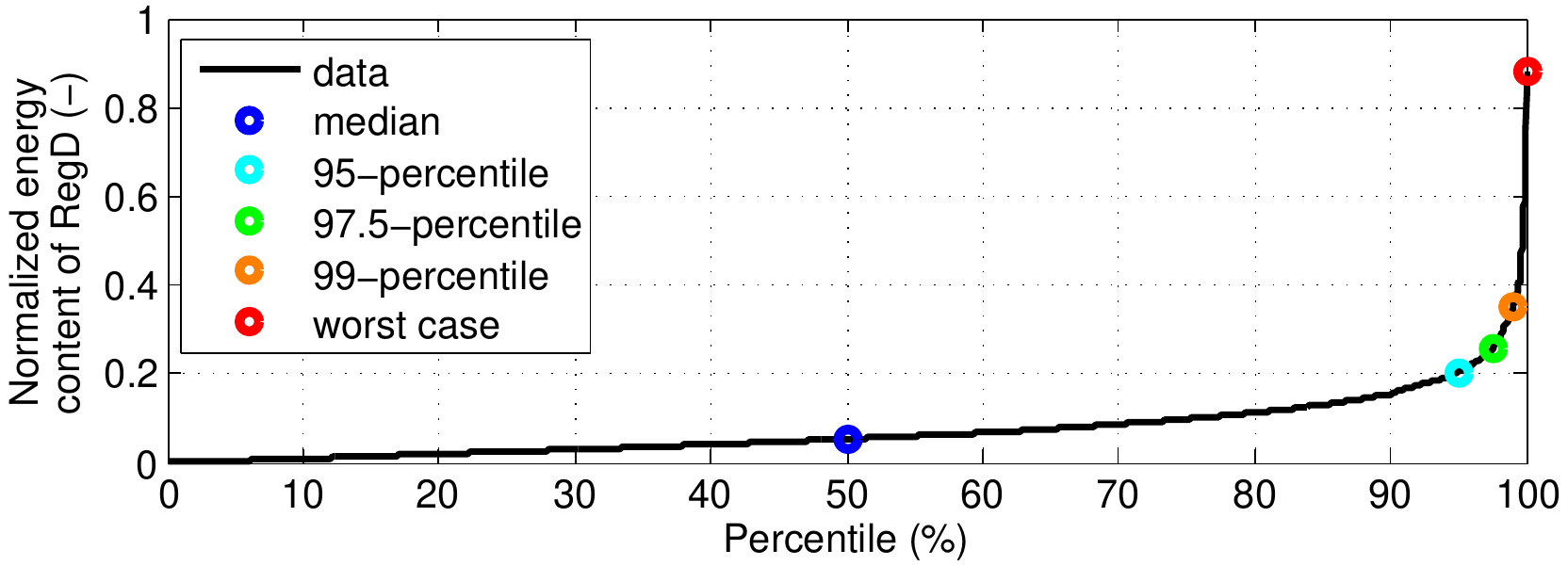}
\caption{Cumulative distribution of RegD signal's energy content in $15$ minutes.} \label{fig:regD_percentiles}
\end{figure}

The energy limit $w_\text{lim}$ of the regulation signal is the worst case normalized reserve request. To identify this limit, we analyze 2-month historical data of the RegD signal from the \ac{PJM} (December $2012$ to January $2013$). Fig.~\ref{fig:regD_percentiles} shows the cumulative distribution of the signal's energy content over $15$ minute intervals, as well as the actual worst case, the median, $95\%$, $97.5\%$, and $99\%$ percentiles. Since the actual worst case of $w_\text{lim}=0.88$ would lead to very conservative solutions, we define the worst case as the $97.5\%$ percentile of the distribution ($w_\text{lim}=0.25$). Therefore, the requested reserve energy by RegD over $15$ minutes will be less than $25\%$ of the reserve capacity with probability $97.5\%$.

The worst case reserve request along the prediction horizon is obtained in \eqref{StateConLv1} by taking the worst case for each time step independently \cite{VrettosIFAC2014}. Although this is a conservative approach, it has the advantage of building additional robustness to modeling and forecast uncertainties.

\subsection{Modeling Reserve Product Constraints} \label{model_res_products}
Problem \eqref{Lv1rca} allows us to select different reserve capacities for each $15$-min time slot, as well as different $R_{\text{u},k}$ and $R_{\text{d},k}$ for the same time slot. However, many markets have requirements on the structure of the reserve product, in particular reserve blocks with minimum duration and/or symmetric reserves.

\emph{Reserve blocks with minimum duration} of $T_\text{res} \in \Nset$ time steps can be modeled by adding in \eqref{Lv1rca} the constraints $R_{\text{u},k}=R_{\text{u},k+j}$ and $R_{\text{d},k}=R_{\text{d},k+j}$ $\forall k=nT_\text{res}+1,~ \forall j \in \{1,\dots T_\text{res}\}$, where $n \in \Nset$ and $n\leq (N_1-1)/T_\text{res}$. We select $T_\text{res}=4$ to require the reserve capacities to be constant over periods of $1$ hour.

\emph{Symmetric reserve capacities} can be enforced by introducing in \eqref{Lv1rca} the constraint $R_{\text{u},k}=R_{\text{d},k}~\forall k$, and are expected to reduce the amount of reserves due to the nonlinear flow-to-power fan model. In addition, symmetric reserves and/or reserve blocks with minimum duration increase the complexity because they are nonlinear equality constraints on $r_{\text{u},k}$ and $r_{\text{d},k}$.

Even if symmetric electric reserve capacities are not required from the resources, we chose to impose symmetry in the thermal domain to limit the impact of offering reserves on the room temperature. If the thermal reserves are symmetric, the electric reserves will be asymmetric due to the nonlinear fan curve and \eqref{Ru_def}, \eqref{Rd_def}. In this case, a single variable $r_k$ can replace $r_{\text{u},k}$ and $r_{\text{d},k}$, which is expected to reduce the computation time. We term this type of reserve offering as \emph{``asymmetric''} operation.

\section{Level $2$: Room Climate Controller} \label{CtrlDesign_lv2}
\subsection{MPC Formulation}
The level $2$ controller determines the air mass flow rate setpoint $u_k$ with the robust MPC formulation
\begin{subequations} \label{Lv2rc}
\begin{align}
	\min_{u_k,r_{\text{u},k},r_{\text{d},k}} &\sum\nolimits_{k=0}^{N_2-1} c_k f(u_k) \label{OFLv2rc}\\
     & u_{\text{min},2} \leq u_k-r_{\text{d},k},~~u_k+r_{\text{u},k} \leq u_{\text{max},2},~\forall k \label{InputConLv2rc} \\
     &R_{\text{u},k}^* = f(u_k)-f(u_k-r_{\text{d},k}), ~\forall k \label{ReserveConLv2rcUp}\\
     &R_{\text{d},k}^* = f(u_k+r_{\text{u},k})-f(u_k), ~\forall k \label{ReserveConLv2rcDown}\\
     &\eqref{StateConLv1rc}, \eqref{DynamicsLv1rcaUp}, \eqref{DynamicsLv1rcaDown}~.\nonumber
\end{align}
\end{subequations}
Problem \eqref{Lv2rc} is similar to \eqref{Lv1rca} with the main differences being (i) the electric reserve capacities $R_{\text{u},k}^*$ and $R_{\text{d},k}^*$ are fixed from level $1$, and (ii) the only objective is to minimize energy cost.

The MPC selects $u_k$, $r_{\text{u},k}$ and $r_{\text{d},k}$ such that the electric reserves $R_{\text{u},k}^*$ and $R_{\text{d},k}^*$ can be provided according to constraints \eqref{ReserveConLv2rcUp} and \eqref{ReserveConLv2rcDown}. Weather forecasts are used in \eqref{Lv2rc} and are updated at every time step. The comfort constraints of \eqref{Lv1rca} and \eqref{Lv2rc} are modeled as soft constraints with high penalties to avoid infeasibility due to plant-model mismatch or forecast errors.

The upper and lower bounds on the air flow rate $u_{\text{min},2}$ and $u_{\text{max},2}$ of the MPC are less tight than those of the reserve scheduler to facilitate meeting the comfort zone constraints
\begin{align}\label{mdot_limits_lv2}
u_{\text{min},2} = h\left(10\%\right), \quad u_{\text{max},2} = h\left(90\%\right)~.
\end{align}
The selected bounds correspond to the minimum acceptable fan speed values suggested by the building manager of FLEXLAB.

\subsection{Kalman Filter}
Since $T_{\text{m},k}$ is not measured directly and the measurement of $T_{\text{r},k}$ is noisy, we use an extended Kalman filter to obtain a state estimate $\hat{x}_k$ for the MPC and the reserve scheduler.

Assuming additive process and measurement noise, the a priori error covariance $P_{\text{e},k}^-$, a posteriori error covariance $P_{\text{e},k}$ and Kalman gain $K_k$ are given by \cite{welch2006kalman}
\begin{align}
P_{\text{e},k}^- &= F_k P_{\text{e},k} F_k^\top + Q \\
K_k &= P_{\text{e},k}^- H_k^\top (H_k P_{\text{e},k}^- H_k^\top + R)^{-1}\\
P_{\text{e},k} &= (I-K_kH_k) P_{\text{e},k}^-~,
\end{align}
where $F_k$ is the Jacobian matrix of system dynamics and $H_k$ is the Jacobian matrix of the output $y_k = C x_k$; $Q=[0.4~0;~0.4~0]$ and $R=0.1$ are the process and measurement noise covariance matrices, respectively;\footnote{We set $R=0.1$ based on the accuracy of the temperature sensors. Based on the building model's RMSE (equal to $0.42^\circ$C from Table~\ref{tab:thermal_model_cmp}), an initial estimate of the diagonal entries of $Q$ is $0.42^2=0.1764$. We chose the larger value $0.4$ because the model's out-of-sample RMSE will be higher than $0.42$.} and $I$ is the identity matrix. Let $\phi_x$ denote the partial derivative of bilinear dynamics and $\varphi_x$ denote the partial derivative of the output equation, both with respect to the state $x_k$. The matrices $F_k$ and $H_k$ are calculated with

\begin{align}
F_k &= \phi_x(\hat{x}_{k-1},u_{k-1}) = A + B_{xu} \hat{x}_{k-1} u_{k-1}\\
H_k &= \varphi_x(\hat{x}_{k-1},u_{k-1}) = C~.
\end{align}

\section{Level $3$: Frequency Regulation Controller} \label{CtrlDesign_lv3}
Level $3$ controls the fan speed (input of the fan controller) such that the fan power tracks the frequency regulation signal. Our approach is different from \cite{MacDonald2014PJM} that used the frequency of the \ac{VFD} as a control variable, and from \cite{LinTSG2014} where a fan speed command was superimposed on the output of the fan controller.

There are four requirements for the frequency regulation controller: fast response, minimal computation effort, accuracy and stability. For this purpose, we developed a novel switched controller with two loops: (i) Ctrl1: a model-based, feedforward controller, and (ii) Ctrl2: a model-free, feedback \ac{PI} controller. The feedforward controller uses the static speed-to-power fan model \eqref{speed-to-power} and it is inherently stable due to the absence of feedback. The PI controller is used to reduce the steady-state error of the feedforward controller, but its stability is not guaranteed and requires gain tuning. The discrete-time implementation of the switched controller is described by Algorithm $1$.

\begin{algorithm}[t]  \label{level3_algorithm}
\begin{algorithmic}[1] 
\small
\State initialize old tracking error $e_{\text{old}} = 0$ and fan speed $N_\text{f}$
\While{experiment is running}
    \State calculate baseline power: $P_\text{s}=f(\dot{m}_\text{s})$
    \State compute reserve: $R\hspace{-0.05cm}=\hspace{-0.05cm}w R_\text{d}~(\text{if}~w\hspace{-0.05cm}>\hspace{-0.05cm}0), R\hspace{-0.05cm}=\hspace{-0.05cm}w R_\text{u}~(\text{if}~w\hspace{-0.05cm}\leq\hspace{-0.05cm}0)$
    \State calculate desired fan power: $P_\text{d} = P_\text{s} + R$
    \Repeat
    \State measure fan power $P_\text{f}$
    \State calculate new tracking error: $e_{\text{new}}=P_\text{d}-P_\text{f}$
    \If{$|e_{\text{new}}| \leq \varepsilon$} \label{errorcheck}
        \State set PI output: $N_\text{f,pi} \hspace{-0.05cm}=\hspace{-0.05cm} N_\text{f} + K_p (e_{\text{new}}\hspace{-0.05cm}-\hspace{-0.05cm}e_{\text{old}}) \hspace{-0.05cm}+\hspace{-0.05cm} K_i \Delta{t} e_\text{new}$
        \State cap fan speed: $N_\text{f}=\text{min}[\text{max}(N_\text{f,pi},N_\text{f,\text{min}}),N_\text{f,\text{max}}]$
        \State set fan speed to $N_\text{f}$
        \State set old tracking error to: $e_{\text{old}} = e_{\text{new}}$
    \Else
        \State set fan speed to: $N_\text{f}=g(P_\text{d})$
        \State set old tracking error to: $e_{\text{old}} = 0$
    \EndIf
    \Until{elapsed time is equal to control loop duration}
\EndWhile
\end{algorithmic}
\caption{Implementation of the switched controller}
\end{algorithm}

Step $3$ of Algorithm $1$ uses the flow-to-power fan model \eqref{flow-to-power} to translate the scheduled flow rate of level $2$ to baseline power consumption. The desired fan power $P_\text{d}$ is computed at step $5$ based on the baseline, the reserve capacity of level $1$ and the regulation signal. The new control error $e_{\text{new}}$ is calculated at step $8$ as the difference between $P_\text{d}$ and the measured fan power $P_\text{f}$. At step $9$ the condition $|e_{\text{new}}| \leq \varepsilon$ is checked to decide whether Ctrl1 or Ctrl2 will be used ($\varepsilon$ is a tolerance that represents the fan model's accuracy). If $|e_{\text{new}}| \leq \varepsilon$ holds, then we activate Ctrl2 (the PI controller's discrete time implementation is given from step $10$ to step $13$). On the other hand, if $|e_{\text{new}}| > \varepsilon$, we activate Ctrl1 and determine the fan speed at step $15$ according to \eqref{speed-to-power}.

After a large power setpoint change, Ctrl1 remains active for as long as $|e_{\text{new}}|$ is larger than $\varepsilon$, whereas the controller switches to Ctrl2 when $|e_{\text{new}}| \leq \varepsilon$. When we switch from Ctrl1 to Ctrl2, we reset the integral error to zero (step $16$ of Algorithm $1$) to avoid large overshoots due to accumulated errors. Furthermore, if the output of Ctrl2 is larger than $90\%$ or smaller than $10\%$, we cap or floor the fan speed to these values.

Due to the nonlinear fan curve, gain scheduling was used in the PI controller. Five operating regions were defined and different proportional ($K_p$) and integral gains ($K_i$) were calculated for each region using the Ziegler-Nichols method \cite{ziegler1942optimum}. We performed step response tests with $K_i=0$ and gradually increased $K_p$ until a critical value with stable and consistent oscillations in fan power. The critical proportional gain and the period of oscillations are used to determine $K_p$ and $K_i$.

The gains obtained with the Ziegler-Nichols method served as an initial guess, whereas the final gains were determined with trial and error and are presented in Table~\ref{tab:PI_tuning}. The $K_p$ gains are lower and the $K_i$ gains are higher than those suggested by the Ziegler-Nichols method because the goal of the PI controller (Ctrl2) is to correct the steady-state error of Ctrl1, but not to recover the system after a large setpoint change.

\begin{table}[t]
\renewcommand{\arraystretch}{1.1}
\caption{Tuned gains of the PI controller}
\centering
\begin{tabular}{c|ccccc}
\hline
Region (kW) & $[0, 0.5)$ & $[0.5, 1)$ & $[1, 1.5)$ & $[1.5, 2)$ & $[2, 2.5)$ \\
$K_p$ (proportional) & 0.004 & 0.004 & 0.004 & 0.0045 & 0.004 \\
$K_i$ (integral) & 0.01 & 0.0035 & 0.003 & 0.0025 & 0.002 \\
\hline
\end{tabular}
\label{tab:PI_tuning}
\end{table}

The proposed switched controller is advantageous in terms of stability and performance compared with the open-loop controller of \cite{MacDonald2014PJM} and the closed-loop controller of \cite{LinTSG2014}. Ctrl1 allows us to track sudden power setpoint changes without the need of high gains in Ctrl2 that would compromise stability.

%% file: conclusion_p1.tex
\section{Conclusion and Outlook} \label{sec_conclusion}
In Part I of this two-part paper, we presented the commercial building test facility FLEXLAB for a frequency regulation demonstration project. We developed and compared different building models for use in a day-ahead reserve scheduler and an MPC for building climate control. Specifically, we presented mathematical reformulations to include the nonlinear fan dynamics in the optimization problems. Furthermore, we proposed a switched controller for frequency regulation that is accurate and inherently stable. In Part II we report extensive experimental results using the developed models and controllers.

%% file: appendix_proofs.tex
\section{Proof of Lemma~\ref{min_max_eq}} \label{proof_lemma2}
\begin{proof}
From the definition of $\Delta{u}_k$ in \eqref{Du_def} and \eqref{res_energy_def} we get
\begin{align} \label{aux_eq1_lemma2}
\underset{w_k}{\min} \left(u_k+\Delta{u}_k \right) &= \underset{w_k \in [-w_\text{lim},0)}{\min} \left[f^{-1}\left(P_k+w_k R_{\text{u},k}\right)\right].
\end{align}
Due to monotonicity, $\text{argmin}[f^{-1}(P_k+w_k R_{\text{u},k})]$ is equal to
\begin{align}
\text{argmin}(P_k+w_k R_{\text{u},k})=-w_\text{lim}.
\end{align}
Substituting the minimizer $-w_\text{lim}$ in \eqref{aux_eq1_lemma2} we get \eqref{lemma2_p1}. Equation \eqref{lemma2_p3} is a special case of \eqref{lemma2_p1} derived as
\begin{align} \label{aux_eq2_lemma2}
&f^{-1}\left(P_k-w_\text{lim} R_{\text{u},k}\right) = f^{-1}\left[f(u_k)-R_{\text{u},k}\right]= \nonumber\\ &f^{-1}\left[f(u_k)-f(u_k)+f(u_k-r_{\text{d},k})\right] = u_k-r_{\text{d},k}~,
\end{align}
where \eqref{Ru_def} is used. The maximization case for \eqref{lemma2_p2} and \eqref{lemma2_p4} can be proved analogously but the proof is omitted for brevity.
\end{proof}

\section{Proof of Proposition~\ref{reformulation_proposition}} \label{proof_propo1}
\begin{proof}
It is sufficient to show that the input constraints \eqref{InputConLv1rc} are equivalent to \eqref{InputConLv1}, and that the set of state constraints \eqref{DynamicsLv1rcUp}, \eqref{DynamicsLv1rcDown} and \eqref{StateConLv1rc} is equivalent to the set of constraints \eqref{DynamicsLv1} and \eqref{StateConLv1}. The equivalence of input constraints follows directly from \eqref{lemma2_p3} and \eqref{lemma2_p4}. For the state constraints we first write \eqref{DynamicsLv1} as
\begin{align}
x_{k+1} \hspace{-0.05cm}=\hspace{-0.05cm} A x_k \hspace{-0.05cm}+\hspace{-0.05cm} (B_u T_{\text{s}} \hspace{-0.05cm}+\hspace{-0.05cm} B_{xu} x_k)\hspace{-0.05cm}\cdot\hspace{-0.05cm}(u_k\hspace{-0.05cm}+\hspace{-0.05cm}\Delta{u}_k) \hspace{-0.05cm}+\hspace{-0.05cm} B_v v_k.
\end{align}
Constraint \eqref{StateConLv1} is applied only to the first state $T_{\text{r},k} = Cx_k$ of the state vector $x_k$, where $C=[1~~0]$ is the output matrix. Observing that $CB_u T_{\text{s}} =b T_{\text{s}}$ and $CB_{xu}x_k = -b T_{\text{r},k}$, constraint \eqref{StateConLv1} can be written as
\begin{align}
\underset{w_k}{\min}\hspace{-0.05cm}\left[CAx_k \hspace{-0.05cm}+\hspace{-0.05cm} b (T_{\text{s}}\hspace{-0.05cm}-\hspace{-0.05cm}T_{\text{r},k}) \hspace{-0.06cm}\cdot\hspace{-0.06cm} (u_k\hspace{-0.05cm}+\hspace{-0.05cm}\Delta{u}_k) \hspace{-0.05cm}+\hspace{-0.05cm} CB_v v_k \right] &\hspace{-0.05cm}\geq\hspace{-0.05cm} x_{\text{min},k} \label{prop1_p1} \\
\underset{w_k}{\max}\hspace{-0.05cm}\left[CAx_k \hspace{-0.05cm}+\hspace{-0.05cm} b (T_{\text{s}}\hspace{-0.05cm}-\hspace{-0.05cm}T_{\text{r},k}) \hspace{-0.06cm}\cdot\hspace{-0.06cm} (u_k\hspace{-0.05cm}+\hspace{-0.05cm}\Delta{u}_k) \hspace{-0.05cm}+\hspace{-0.05cm} CB_v v_k \right] &\hspace{-0.05cm}\leq\hspace{-0.05cm} x_{\text{max},k}. \label{prop1_p2}
\end{align}
Due to assumption~\ref{SAT_assumption}, we have $b (T_{\text{s}}-T_{\text{r},k})\leq0$ and so the left hand side of \eqref{prop1_p1} is minimized when $u_k+\Delta{u}_k$ is maximized. From \eqref{lemma2_p2}, this is achieved when $u_k+\Delta{u}_k = f^{-1}\left(P_k+w_\text{lim} R_{\text{d},k}\right)$ holds. Let $\underline{x}_{k}$ denote the minimum state trajectory at time step $k$. The time evolution of $\underline{x}_{k}$ is therefore obtained by $\underline{x}_{k+1} = CA \underline{x}_{k} + b (T_{\text{s}}-C\underline{x}_{k}) \cdot f^{-1}\left(P_k+w_\text{lim} R_{\text{d},k}\right) + CB_v v_k$, which is essentially \eqref{DynamicsLv1rcDown}. Thus, if $x_{\text{min},k} \leq \underline{x}_{k}$ holds, then the left hand side of \eqref{StateConLv1} also holds because $x_k \leq \underline{x}_{k}$~$\forall w_k$.

Analogous arguments can be used to show that \eqref{prop1_p2} results in \eqref{DynamicsLv1rcUp} and $\overline{x}_{k} \leq x_{\text{max},k}$, which are equivalent to the right hand side of \eqref{StateConLv1}, but this is omitted for brevity.
\end{proof}

\section{Proof of Proposition~\ref{bound_proposition}} \label{proof_propo2}
\begin{proof}
Using definition \eqref{Ru_def} and the convexity of $f$ we get
\begin{align} \label{aux1_propo3}
P_k-w_\text{lim} R_{\text{u},k} &= f(u_k)-w_\text{lim}\left[f(u_k)-f(u_k-r_{\text{d},k})\right]\nonumber\\
&=(1-w_\text{lim}) f(u_k) + w_\text{lim} f(u_k-r_{\text{d},k}) \nonumber\\
&\geq f\left[ (1-w_\text{lim}) u_k + w_\text{lim} (u_k-r_{\text{d},k}) \right]~.
\end{align}
Using the monotonicity of $f$ in \eqref{aux1_propo3} we get
\begin{align} \label{aux2_propo3}
f^{-1}\left(P_k-w_\text{lim} R_{\text{u},k}\right) &\geq (1-w_\text{lim}) u_k + w_\text{lim} (u_k-r_{\text{d},k})\nonumber\\
&= u_k-w_\text{lim} r_{\text{d},k}~.
\end{align}
$\overline{x}_k^* \geq \overline{x}_k^\star$ is now obtained by combining \eqref{aux2_propo3}, \eqref{DynamicsLv1rcUp}, \eqref{DynamicsLv1rcaUp}, and using the same arguments related to $b (T_{\text{s}}-T_{\text{r},k})\leq0$ as in the proof of Proposition~\ref{reformulation_proposition}. Similarly, one can show that $u_k+w_\text{lim} r_{\text{u},k} \leq f^{-1}\left(P_k+w_\text{lim} R_{\text{d},k}\right)$ holds and so that $\underline{x}_k^* \geq \underline{x}_k^\star$ also holds. Fig.~\ref{fig:Fan_model_fit_explain} provides a graphical interpretation of \eqref{aux2_propo3}.
\end{proof}